\documentclass[DIV=calc, paper=letter, fontsize=11pt]{scrartcl}
\usepackage[]{units}
\usepackage{lipsum} 
\usepackage[english]{babel} 
\usepackage[protrusion=true,expansion=true]{microtype} 
\usepackage{amsmath,amssymb,amsfonts,amsthm} 
\usepackage[svgnames, table]{xcolor} 
\usepackage[hang, small,labelfont=bf,up,textfont=it,up]{caption} 
\usepackage{booktabs} 
\usepackage{fix-cm}	 
\usepackage[]{units}

\usepackage{bm}
\usepackage{tabularx} 

\usepackage[colorlinks=true, linkcolor=blue, citecolor=blue, filecolor=blue, runcolor=blue, urlcolor = blue]{hyperref}

\usepackage{fullpage}
\usepackage{sectsty} 

\usepackage{fancyhdr} 
\usepackage{lastpage} 

\usepackage{nicefrac}

\usepackage{tikz}
\usetikzlibrary{calc,patterns,decorations.pathmorphing,decorations.markings}
\usepackage{pgfplots}

\usepackage{stmaryrd}
\usepackage{multicol}
\usepackage{float}

\usepackage{subcaption}

\newcounter{tempcolnum}
\makeatletter
\newcommand{\multicolinterrupt}[1]{
\setcounter{tempcolnum}{\col@number}
\end{multicols}
#1%
\begin{multicols}{\value{tempcolnum}}
}
\makeatother

\definecolor{issuePJA_color}{rgb}{1.0,0.0,0.0}

\definecolor{commentPJA_color}{rgb}{1.0,0.0,0.8}

\definecolor{issueNAT_color}{rgb}{0.0,0.5,0.0}

\definecolor{commentNAT_color}{rgb}{0.0,0.6,0.3}

\definecolor{commentWP_color}{rgb}{0.0,0.2,0.6}

\newcommand{\lla}{\left\langle}
\newcommand{\rra}{\right\rangle}
\newcommand{\llb}{\llbracket}
\newcommand{\rrb}{\rrbracket}

\definecolor{rev_color}{rgb}{0.6,0.0,0.0}

\usepackage{graphicx}
\graphicspath{{figData/}}
\DeclareGraphicsExtensions{.png}

\usepackage{lettrine} 
\newcommand{\initial}[1]{ 
\lettrine[lines=3,lhang=0.3,nindent=0em]{
\color{DarkGoldenrod}
{\textsf{#1}}}{}}

\usepackage{caption}

\usepackage{graphicx}

\usepackage{titling} 

\newcommand{\HorRule}{\color{DarkGoldenrod} \rule{\linewidth}{1pt}} 

\pretitle{\vspace{-80pt} \begin{flushleft} \HorRule \fontsize{17}{17} \color{DarkRed} \selectfont} 

\title{Stochastic Discontinuous Galerkin Methods (SDGM) Based on Fluctuation-Dissipation Balance}
\posttitle{\par\end{flushleft}} 
\preauthor{\vspace{-5pt} \begin{flushleft}\fontsize{14}{14} \large  \color{DarkRed}} 
\author{W.\ Pazner$^{*}$, N.\ Trask$^{**}$, and P.\ J.\ Atzberger$^{\dagger}$ }
\postauthor{\fontsize{12}{12} \small \color{Black} 
$*$ Department of Mathematics, University of California, Berkeley; $**$ Sandia National Laboratories;
$\dagger$ Department of Mathematics, University of California, Santa Barbara; atzberg@gmail.com; http://atzberger.org/ 
\vskip 1.5em
\lfoot{\footnotesize * Work supported for P.J.A., N.T. by DOE ASCR CM4 DE-SC0009254 and NSF-MSPRF program, \\ and P.J.A. by NSF Grant DMS - 1616353.  W.P. by Department of Defense, National Defense Science and Engineering Graduate. Fellowship Program and Natural Sciences and Engineering Research Council of Canada. Sandia National Laboratories is a multi-mission laboratory managed and operated by National Technology and Engineering Solutions of Sandia, LLC., a wholly owned subsidiary of Honeywell International, Inc., for the U.S. Department of Energy's National Nuclear Security Administration under contract DE-NA0003525.}
\end{flushleft} \vspace{-18pt} \HorRule} 

\newcommand{\mb}[1]{\mathbf{#1}}
\newcommand{\bs}[1]{\boldsymbol{#1}}

\DeclareGraphicsExtensions{.png}

\begin{document}

\maketitle 

\thispagestyle{fancy} 



\initial{W}\textbf{e introduce a general framework for approximating parabolic Stochastic Partial Differential Equations (SPDEs) based on fluctuation-dissipation balance.  Using this approach we formulate Stochastic Discontinuous Galerkin Methods (SDGM).  We show how methods with linear-time computational complexity can be developed for handling domains with general geometry and generating stochastic terms, handling both Dirichlet and Neumann boundary conditions.  We demonstrate our approach on example systems and contrast with alternative approaches using direct stochastic discretizations based on random fluxes.  We show how our \textit{Fluctuation-Dissipation Discretizations} (FDD) framework allows for compensating for differences in dissipative properties of discrete numerical operators relative to their continuum counter-parts.  This allows us to handle general heterogeneous discretizations, accurately capturing statistical relations.  Our FDD framework provides a general approach for formulating SDGM discretizations and other numerical methods for robust approximation of stochastic differential equations. }

\vspace{12pt}

\setlength{\parindent}{5ex}

\section{Introduction}
\label{sec:Intro}

Stochastic Partial Differential Equations (SPDEs)~\cite{Hairer2009,bookSPDEsDaPrato2014,bookSPDEsLiu2015,Karniadakis2017} arise in many settings including statistical inference~\cite{Cialenco2018,Lototsky2009,Koski1985}, reduced descriptions of dynamical systems~\cite{AtzbergerTabak2015,Fox1970,Stinis2018,Zwanzig1961,Mori1965}, and applications in the natural sciences and engineering~\cite{AtzbergerSELM2011,AtzbergerSIB2007,Fox1970,Donev2010,Uma2011,LandauFld1986}.  Applications include formulations of stochastic phase-field equations~\cite{DaPratoStochCahnHilliard1996,Oden2015,Bertini2002,Sancho1998}, stochastic concentration equations~\cite{bookSPDEsDaPrato2014,AtzbergerRD2010,Donev2017}, and fluctuating hydrodynamics~\cite{LandauFld1986,Donev2010,Donev2014,Atzberger2006,AtzbergerSIB2007,Karniadakis2007,BalboaUsabiaga2012} with fluid-structure interactions~\cite{AtzbergerSIB2007,AtzbergerFoundationSIB2008, AtzbergerSELM2011,AtzbergerLAMMPS2016}.  The development of effective computational methods for SPDEs poses significant challenges often not encountered in the deterministic setting.  The solutions of SPDEs take on the form of a measure on a space of functions or more often on distributions or generalized functions (separable Banach spaces)~\cite{Lieb2001,Hairer2009,bookSPDEsDaPrato2014}.  Stochastic equations are often driven by Gaussian noise in time and space with statistical structures ranging from $\delta$-correlations (which are often not well-posed)~\cite{bookSPDEsDaPrato2014,Karniadakis2017} to colored noise with prescribed spatial-temporal correlation functions~\cite{bookSPDEsLiu2015,AtzbergerSELM2011}.  Such stochastic fields, when they do exist, are often non-differentiable in time, or do not have well-defined pointwise values, which requires interpretations in terms of a measure on an appropriate space of functions or distributions~\cite{bookSPDEsLiu2015,Hairer2009}.  These issues are further compounded in numerical methods when attempting to formulate approximations of the stochastic terms in the presence of truncation errors that arise from discretizations of the underlying differential operators~\cite{AtzbergerSIB2007,AtzbergerRD2010}.

Many applications require preservation of inherent statistical structures to accurately capture stationary quantities or thermodynamic properties.  We shall develop numerical methods to do this by driving them stochastically in a manner that compensates for discretization errors.  We develop our approach utilizing considerations of fluctuation-dissipation balance.  The fluctuation-dissipation principle is a well-established result in equilibrium statistical mechanics having its historic origins in the investigation of noise in electrical circuits (Johnson-Nyquist noise)~\cite{Johnson1928}.  For circuits, it was observed that the amplitude of fluctuations is closely related to electrical impedance~\cite{Nyquist1928}.  These observations subsequently led to the formulation of general principles for physical systems relating the amplitude of fluctuations to the linear responses to perturbations~\cite{Callen1951,Onsager1931,Reichl1998}.  We abstract these notions to the numerical setting by considering how fluctuations should relate to numerically captured dissipative relaxations of the system.  This allows us to obtain discretization-dependent relations, which in turn can be used to develop strategies for the prescription of stochastic driving fields in numerical methods that preserve inherent statistical structures.  These relations allow for taking into account the discretization errors and resulting particularities of the dissipative properties of the utilized discrete numerical operators.  We have successfully used related ideas in our prior works when formulating stochastic numerical methods in~\cite{AtzbergerSIB2007, AtzbergerRD2010, AtzbergerFEM2014, AtzbergerShear2013}. 

\lfoot{\footnotesize }

In this work, we abstract these ideas to formulate a general mathematical framework for stochastic numerical methods for SPDEs which we refer to as \textit{Fluctuation Dissipation Discretizations} (FDD).  For spatial-temporal numerical discretizations, we derive conditions relating the covariances of the stochastic driving terms to the corresponding covariances of the fluctuations in the statistical steady-state.  We show how this FDD framework can be used to develop general strategies for prescribing spatial or temporal discretizations for stochastic equations that control discretization artifacts or attain discrete statistical structures amenable to efficient computational methods.  

In practice, a central challenge is efficiently generating the prescribed stochastic driving fields.  For finite differences with uniform and staggered mesh discretizations, the FDD framework shows covariances can be factored in terms of divergence and gradient mesh operators recovering methods for efficient stochastic driving terms related to random fluxes~\cite{AtzbergerSIB2007,balboa2012staggered,AtzbergerWangSELMChannel2018}.  For non-uniform discretizations, as arise in finite volume and finite element methods, FDD gives much more complicated covariances and work has been done to address this by employing factorizations and multigrid techniques to obtain efficient Gibbs samplers to generate stochastic driving fields~\cite{AtzbergerRD2010,AtzbergerFEM2014}.   Here, we take a different approach for non-uniform discretizations and general geometries avoiding the need for multigrid-based Gibbs samplers.  We show how a judicious choice of spatial discretization can be used to achieve efficient stochastic methods.  We develop a family of Stochastic Discontinuous Galerkin Methods (SDGM) based on FDD, which are applicable to general unstructured curvilinear grids and boundary conditions.  A key feature of discretization a natural block-structure of the prescribed FDD covariances, which can be readily factored.  This allows us to develop stochastic methods that have linear-time computational complexity $O(N)$ under $h$-refinement, which are also amenable also to parallelization.

We remark that our approaches based on fluctuation-dissipation balance are similar in spirit to other recent developments in numerical analysis following a long trend motivated by mimicking key inherent features of the mathematical problem being approximated~\cite{arnold2007compatible,Arnold2010}.  This includes mimetic methods (MM)~\cite{Hyman2002}, discrete exterior calculus (DEC)~\cite{DesbrunHirani2003}, and finite element exterior calculus (FEEC)~\cite{Arnold2010,HolstGilette2017}. Each of these, in their own way, aim to preserve in the numerics inherent geometric features, relations in mechanics, or other properties important to the application domain. In the deterministic setting, this is often observed to be essential to obtaining efficient convergent methods, especially for finite element methods~\cite{arnold2007compatible,Arnold2010,HolstGilette2017}.  Here, our introduced FDD framework provides a way to preserve inherent statistical structures important when approximating stochastic systems.

We organize our paper by first discussing some background on stochastic partial differential equations (SPDEs) in Section~\ref{sec:spde}.  We then formulate our \textit{Fluctuation Dissipation Discretizations} (FDD) framework and show how the FDD framework may be used to develop general strategies for prescribing spatial or temporal discretizations for approximating stochastic equations in Section~\ref{sec:fluct_dissip}.  We then show how these ideas can be used to formulate Stochastic Discontinuous Galerkin Methods (SDGM) in Section~\ref{sec:stoch_galerkin}.  We develop computational methods having linear-time computational complexity for generating stochastic driving fields obeying our FDD framework in Section~\ref{sec:discont_galerkin} and~\ref{sec:stoch_driving_fields}.  We then make comparisons between our FDD-SDGM methods with the alternative approach of directly discretizing the stochastic terms using random fluxes in Section~\ref{sec:results}.  These results highlight the utility of our FDD approach in providing a means to control discretizations artifacts and spurious long-range correlations.  We then present results showing how our FDD-SDGM methods can handle general geometries and periodic, Neumann, and Dirichlet boundary conditions in Sections~\ref{sec:results_periodic}--\ref{sec:results_dirichlet}.  The results show how our FDD approach can be used to devise general strategies for mitigating artifacts from the underlying spatial-temporal discretizations to develop robust numerical methods for stochastic differential equations.

\section{Stochastic Partial Differential Equations}
\label{sec:spde}

We consider parabolic stochastic partial differential equations (SPDEs) of the form 
\begin{eqnarray}
\label{equ:spde_general}
\frac{\partial \mb{u}}{dt} = \mathcal{L} \mb{u} + \mb{b} + \mb{g},
\end{eqnarray}
where $\mathcal{L}$ is a uniform elliptic operator, $\mb{b}$ source term that is square-integrable, and $\mb{g}$ a stochastic force~\cite{bookSPDEsDaPrato2014,bookSPDEsLiu2015}.  We take $\mb{g}$ to be a Gaussian process that is $\delta$-correlated in time with mean zero $\langle \mb{g} \rangle = 0$ and prescribed covariance structure $\mathcal{G}$
\begin{eqnarray}
\langle \mb{g}(\mb{x},t) \mb{g}(\mb{y},s) \rangle 
 =  \mathcal{G}(\mb{x},\mb{y})\delta(t - s).
\end{eqnarray}
As a specific example, we shall consider the stochastic diffusion equation capturing linearized concentration fluctuations~\cite{AtzbergerRD2010}, which can be expressed in this form as
\begin{eqnarray}
\label{equ:stoch_diffusion_main}
\frac{\partial {u}}{dt} = D\Delta {u} + {b} + {g},
\end{eqnarray}
where  
\begin{eqnarray}
\langle {g}(\mb{x},t) {g}(\mb{y},s) \rangle 
 =  -2 \bar{u} D\Delta_{\mb{x}} \delta(\mb{x}-\mb{y}) \delta(t - s).
\end{eqnarray}
This corresponds to $\mathcal{L} = D\Delta$ and spatial correlations $\mathcal{G}(\mb{x},\mb{y}) = -2 \bar{u} D\Delta_{\mb{x}} \delta(\mb{x}-\mb{y})$.  Here, $D$ denotes the diffusivity, and $\bar{u}$ the reference concentration around which the system was linearized.  This has $\delta$-correlations in time and correlations $\Delta_{\mb{x}} \delta$ in space~\cite{AtzbergerRD2010}.  The statistical steady-state solution has equilibrium fluctuations $\mathcal{C}(\mb{x},\mb{y}) = \langle {u}(\mb{x}){u}(\mb{y}) \rangle = \bar{u} \delta(\mb{x} - \mb{y})$.  As a consequence, the solutions of such stochastic differential equations are irregular, being non-differentiable in time and not having point-wise values but rather solutions in the sense of a measure on a space of distributions~\cite{Hairer2009,bookSPDEsDaPrato2014,bookSPDEsLiu2015,Oksendal2000}.  This poses significant challenges both for analysis and numerical approximations.  Many approaches for SPDEs rely on spectral approximations of the solutions~\cite{StuartHairerSpectralSPDEs2016,Karniadakis2017,bookSPDEsDaPrato2014}.  We take a different approach, using instead finite differences and finite elements to discretize the differential operators, temporal dynamics, and stochastic driving fields~\cite{AtzbergerSIB2007,AtzbergerRD2010,
AtzbergerFEM2014,AtzbergerShear2013,Du2002,Allen1998,Donev2014}.

For the linear stochastic differential equations~\ref{equ:spde_general}, the solutions $\mb{u}$ are Gaussian processes and have significant statistical structures that we can utilize in constructing our numerical approximations.  The Gaussianity ensures that the process is determined if we know the mean values and spatial-temporal correlation functions of the process~\cite{Hairer2009,bookSPDEsLiu2015}.  The solution can be expressed formally as
\begin{equation}
\label{equ:formal_sol_u}
u(\mb{x},t) 
 = \exp\left(-t\mathcal{L}\right)u(0) + \int_0^{t} \exp\left(-(t-s)\mathcal{L}\right) b(s) ds + \int_0^{t} \exp\left(-(t-s)\mathcal{L}\right) \mathcal{G}(\mb{x},\mb{y}) d\mathcal{W}_s. \hspace{0.1cm}
\end{equation}
Here, $\exp\left(-t\mathcal{L}\right)$ denotes the semigroup of solution operators for the parabolic PDE with operator $\mathcal{L}$ and $\mathcal{W}_s = \mathcal{W}(\mb{x},s)$ is a formal space-time Brownian motion with integral given the Ito interpretation~\cite{Hairer2009,bookSPDEsLiu2015,Oksendal2000}.  The mean of this process is given formally by
\begin{eqnarray}
\langle u(\mb{x},t) \rangle 
 =  \exp\left(-t\mathcal{L}\right)\langle u(0) \rangle  + \int_0^{t} \exp\left(-(t-s)\mathcal{L}\right) \langle b(s) \rangle ds.
\end{eqnarray}
When $b = 0$, the stationary distribution has $\langle u(0) \rangle = 0$.  The spatial-temporal correlation are then
\begin{eqnarray}
\langle u(\mb{x},t) u(\mb{y},s) \rangle 
 =  \exp\left(-(t-s)\mathcal{L}_{\mb{x}}\right)\mathcal{C}(\mb{x},\mb{y}),
\end{eqnarray}
where $t \geq s$, and $\mathcal{C}(\mb{x},\mb{y}) = \langle u(\mb{x},0) u(\mb{y},0) \rangle$ is the steady-state covariance of fluctuations of $u$. This can be generalized to non-zero $b$ by subtracting off the mean behavior to consider $\tilde{u} = u - \bar{u}$ with $\bar{u}(t) = \langle u(t) \rangle$.  Since the process is Gaussian we can formally express the probability measure in terms of a formal density as 
\begin{eqnarray}
\label{equ:spde_formal_density}
\rho[u] = (1/\mathcal{Z})\exp\left(-\frac{1}{2}\mb{u}^T \mathcal{C}^{-1}\mb{u}\right),
\end{eqnarray}
where $\mathcal{Z}$ is the normalization factor.  This should be interpreted as defining the statistics of finite dimensional marginals~\cite{bookSPDEsLiu2015}.  

\subsection{Spatial Discretization}
\label{sec:spde_discr_spatial}

We approximate solutions of the stochastic partial differential equation~\ref{equ:spde_general} using a discrete stochastic dynamical system of the form
\begin{eqnarray}
\label{eq:governingEq_discr}
\frac{d{Z}_t}{dt} = L {Z}_t + F({Z}_t)+ {F}_{noise}.
\end{eqnarray}
We denote by $Z_t \in \mathbb{R}^n$ the state of the system with $Z_t \sim u(t)$, by $F({Z}_t)$ the system forcing with $F({Z}_t) \sim b$, and by ${F}_{noise} = Q \frac{dW_t}{dt}$ the stochastic driving force with ${F}_{noise} \sim \mb{g}$.  We take ${F}_{noise}$ to be a Gaussian process with mean zero and a yet-to-be-determined covariance 
\begin{eqnarray}
\langle {F}_{noise}(s) {F}_{noise}^T(t) \rangle 
= G \delta(s - t),
\end{eqnarray}
where $G = QQ^T$.  The solution $Z_t$ is a Gaussian process completely determined by its mean and spatial-temporal covariances.  Similar to equations~\ref{equ:formal_sol_u}--\ref{equ:spde_formal_density}, we have that the system has the mean 
\begin{eqnarray}
\langle Z_t \rangle 
 = \exp\left(-t{L}\right)\langle Z_0 \rangle  + \int_0^{t} \exp\left(-(t-s)\mathcal{L}\right) \langle F(Z_s) \rangle ds.
\end{eqnarray}
When $F = 0$ we have $\langle Z_t \rangle = 0$ and the temporal correlations are
\begin{eqnarray}
\label{equ:autocor_discr_spatial}
\langle Z_t Z_s^T \rangle = \exp\left(-(t-s)L\right) C,
\end{eqnarray}
where $t > s$.  For the covariance $C_t = \langle Z_t \left(Z_t\right)^T \rangle$, we assume throughout that we have ergodicity so that $C_t \rightarrow C$ as $t \rightarrow \infty$ with the steady-state covariance $C = \langle Z Z^T \rangle$.  This can be generalized when $f$ is non-zero by subtracting off the mean $\tilde{Z}_t = Z_t - \bar{Z}_t$ with $\bar{Z}_t = \langle Z_t \rangle$.

\subsection{Temporal Discretization}
\label{sec:spde_discr_temporal}

We shall also consider the role of temporal discretization errors.  We consider for illustration the case of an Euler discretization
\begin{eqnarray}
\label{equ:discr_time_equ}
{Z}^{n+1} = {Z}^{n} + \Delta{t}L{Z}^n + \Delta{t}{F}^n + Q \Delta W^n,
\end{eqnarray}
where $\Delta W^n$ is the Brownian increment $W^{n+1} - W^n$. This is a discrete-time process (multi-variate Gaussian for finite number of steps) which is determined by its mean and covariance.  The mean is
\begin{eqnarray}
\langle {Z}^{n} \rangle  = \left(I + \Delta{t}L\right)^n \langle {Z}^{0} \rangle + \sum_{k=0}^{n-1} \left(I + \Delta{t}L\right)^k \Delta{t} \langle {F}^{n - 1 - k} \rangle.
\end{eqnarray}
In the case when $F^{s} = 0$ we have $\langle {Z}^{n} \rangle = 0$.  The covariance is
\begin{eqnarray}
\label{equ:autocor_discr_time}
\langle {Z}^{n} \left({Z}^{m}\right)^T \rangle = \left(I + \Delta{t}L\right)^{n - m} C
\end{eqnarray}
where $n\geq m$, and $C = \langle {Z}{Z}^T \rangle$ is the steady-state covariance. This can be generalized for non-zero $F^s$ by subtracting the mean $\tilde{Z}^n = Z^n - \bar{Z}^n$ with $\bar{Z}^n = \langle Z^n \rangle$.

For these spatial-temporal discretization cases, we show how we can introduce stochastic driving terms $F_{noise}(t)$ or ${F}_{noise}^n$ that compensate in numerical methods for spatial truncation errors or those arising from the choice of temporal approximations and numerical integrator.

\section{Fluctuation-Dissipation Discretizations (FDD) Framework}
\label{sec:fluct_dissip}

We take the general approach for discretizations of the stochastic terms in equation~\ref{equ:spde_general} of using variations of fluctuation-dissipation balance adapted to the setting of numerical discretizations in space or in time. We express equation~\ref{eq:governingEq_discr} representing a semi-discretization with continuous time as 
\begin{eqnarray}
\label{equ:sde_dZt}
dZ_t = LZ_t dt + QdW_t.
\end{eqnarray}
We give this the interpretation of an Ito Stochastic Process~\cite{Oksendal2000}.  Letting $C_t = \langle Z_t Z_t^T \rangle$, we have $dC_t = \langle dZ_t Z_t^T \rangle + \langle Z_t dZ_t^T \rangle$.  We have from Ito's Lemma~\cite{Oksendal2000} and equation~\ref{equ:sde_dZt} that
\begin{eqnarray}
dC_t = LC_t + C_t^TL^T + G
\end{eqnarray}
where $G = QQ^T$.  This gives in the statistical steady-state with $dC_t \rightarrow 0$ as $t \rightarrow \infty$ the relationship
\begin{eqnarray}
\label{equ:fluct_dissp_discr_space}
G = - LC - (LC)^T.
\end{eqnarray}
This establishes a fluctuation-dissipation relationship for the discrete system given in equations~\ref{eq:governingEq_discr} and~\ref{equ:sde_dZt}.  Throughout, we assume that $L$ is negative-definite and that the stochastic driving fields yield ergodic behaviors for the stochastic process.  

In the case that $LC = (LC)^T$, we have
\begin{eqnarray}
\label{equ:C_relate_G_spatial}
C = -\frac{1}{2}L^{-1}G.
\end{eqnarray}
If we make a specific choice for $G$ for the stochastic driving terms for the system, this establishes what equilibrium fluctuations $C$ will be obtained.

We can also obtain analogous relations when time has been discretized as in equation~\ref{equ:discr_time_equ}.  In this case, we have for $C^{n+1} = \langle Z^{n+1} \left(Z^{n+1}\right)^T \rangle$ that
\begin{eqnarray}
C^{n+1} = C^{n} + \Delta{t}LC^n + \Delta{t}C^nL^T + \Delta{t}G + \Delta{t}^2LC^nL^T.
\end{eqnarray}
This was obtained by substituting for $Z^{n+1}$ using equation~\ref{equ:discr_time_equ} with $F^n = 0$.  In the statistical steady-state, we have $C^{n} - C^{n+1} \rightarrow 0$ as $n \rightarrow \infty$.  This gives the fluctuation-dissipation relation taking temporal discretization into account
\begin{eqnarray}
\label{equ:fluct_dissp_discr_time}
G = -LC - (LC)^T - \Delta{t}LCL^T.
\end{eqnarray}
We see the temporal discretization results in a higher-order correction for the stochastic driving fields involving the time-step size to compensate for the additional truncation errors of the numerical methods.

In the case that $LC = (LC)^T$, we have the further useful relation
\begin{eqnarray}
\label{equ:C_relate_G_time}
C = -\frac{1}{2}L^{-1}\left(I + \frac{1}{2}\Delta{t}L^T\right)^{-1}G.
\end{eqnarray}
For a specific choice for $G$ for the stochastic driving terms for the system, this establishes what equilibrium fluctuations $C$ will be obtained.  We see how the temporal discretization augments the obtained equilibrium fluctuations relative to equation~\ref{equ:C_relate_G_spatial}.

The fluctuation-dissipation relations in equation~\ref{equ:fluct_dissp_discr_space} and~\ref{equ:fluct_dissp_discr_time} provide a few strategies for discretizing stochastic equations.  For many stochastic equations arising from applications in physics the equilibrium fluctuations are known from the energy of the system $E[u]$.  From equation~\ref{equ:spde_formal_density}, when $E[u]$ is quadratic in $u$, or when the system is linearized, this yields the equilibrium covariance $\mathcal{C}$.  When discretizing the energy $E[u] \rightarrow E[Z]$ this yields a covariance $C$ for the discretized system.  One strategy is to discretize the stochastic driving terms $\mb{g}$ by using relations in equation~\ref{equ:fluct_dissp_discr_space} or~\ref{equ:fluct_dissp_discr_time} to obtain a covariance $G$ for the stochastic driving terms in the discrete system. This ensures the stochastic driving terms interact with the discrete numerical operator L to yield equilibrium fluctuations with covariance C, even in the presence of discretization errors.

An alternative strategy is to use equation~\ref{equ:C_relate_G_spatial} or~\ref{equ:C_relate_G_time} to formulate $G$ that yields an equilibrium covariance $C$ with acceptable properties such as localization in space with rapid decay to approximate $\delta$-correlations or reduce artifacts at coarse-refined interfaces within structured discretizations~\cite{AtzbergerRD2010}.  Without these considerations, arbitrarily approximating directly $\mb{g}$ can result in ringing phenomena or other artifacts in the numerics~\cite{AtzbergerRD2010}.  We also see from equations~\ref{equ:autocor_discr_spatial} and~\ref{equ:autocor_discr_time} the relations can be used to help ensure accurate autocorrelation functions for the discretized stochastic system.  

We have used related ideas in our prior works~\cite{AtzbergerSIB2007,AtzbergerSELM2011,AtzbergerTabak2015} to derive stochastic discretizations for finite difference methods on uniform meshes, staggered meshes, and spatially adaptive meshes~\cite{AtzbergerSIB2007,AtzbergerWangSELMChannel2018,AtzbergerRD2010} and for finite element methods~\cite{AtzbergerFEM2014}.  In the subsequent sections, we utilize the FDD framework and related ideas to develop Stochastic Discontinuous Galerkin Methods (SDGM) to handle domains with general geometries and Neumann and Dirichlet boundary conditions.

\section{Stochastic Discontinuous Galerkin Methods (SDGM)}
\label{sec:stoch_galerkin}

We develop Stochastic Discontinuous Galerkin Methods (SDGM) by introducing stochastic driving terms for approximating solutions of the SPDE in equation~\ref{equ:spde_general}. The Discontinuous Galerkin (DG) method was first introduced in 1973 by Reed and Hill for solving the neutron transport equation \cite{Reed1973}. Subsequently, Cockburn and Shu extended the DG method
to general systems of nonlinear hyperbolic conservation laws \cite{Cockburn1991,
Cockburn1989a, Cockburn1989b, Cockburn1990, Cockburn1998}. Several extensions have been developed for elliptic and parabolic problems 
\cite{Arnold1982,Bassi1997,Brezzi2000,Cockburn1998a}. These extensions subsequently have been presented in the context of a unified framework in \cite{Arnold2002}.

\subsection{Discontinuous Galerkin (DG) Methods}
\label{sec:discont_galerkin}

The DG method is a finite element method, suitable for use on unstructured
meshes, that makes use of a discontinuous, piecewise polynomial function space.
To begin, we describe the DG discretization of the the deterministic heat
equation on a spatial domain $\Omega \subseteq \mathbb{R}^2$,
\begin{align}
  \label{eq:heat-eqn-1}
  \frac{\partial u}{\partial t} &= \Delta u + f, && \text{in $\Omega$},\\
  \label{eq:heat-eqn-2}
  u &= g_D, && \text{on $\Gamma_D$},\\
  \label{eq:heat-eqn-3}
  \frac{\partial u}{\partial \bm n} &= \bm g_N \cdot \bm n,
    && \text{on $\Gamma_N$},
\end{align}
where Dirichlet and Neumann boundary conditions are imposed on the boundary
$\Gamma_D \cup \Gamma_N = \partial\Omega$.  We discretize the spatial domain
$\Omega \subseteq \mathbb{R}^2$ by introducing a mesh
\begin{equation}
  \mathcal{T}_h = \left\{
    K_i, 1 \leq i \leq n_t, \bigcup_{i=1}^{n_t} K_i = \Omega
  \right\}.
\end{equation}
One of the main advantages of the DG method is the great amount of flexibility
it affords in the construction of the mesh. The mesh can be be entirely
unstructured, non-conforming, and can consist of elements of different shapes.
In this work, we restrict the elements of the mesh $K_i$ to be straight-sided
quadrilaterals. The generalization to curved elements is possible using a
standard isoparametric mapping \cite{Hesthaven2008}. We now fix a polynomial
degree $p \geq 0 $. On each element $K_i$, we define the space of bivariate
polynomials of degree at most $p$ in each variable,
\begin{equation}
  \mathcal{Q}^p(K_i) = \left\{
    v(\bm x) : K_i \to \mathbb{R} :
    v(\bm x) = \sum_{\bm\alpha} c_\alpha \bm x^{\bm\alpha}
  \right\},
\end{equation}
where the sum is taken over all multi-indices $\bm\alpha = (\alpha_1, \alpha_2)$
such that $\alpha_i \leq p$, and the notation $\bm x^{\bm\alpha}$ is used to 
mean $x_1^{\alpha_1} x_2^{\alpha_2}$. We now introduce the discontinuous finite
element space $V_h$ defined by
\begin{equation}
  V_h = \left\{ v_h(\bm x) : \Omega \to \mathbb{R} : 
  v_h |_{K_i} \in \mathcal{Q}^p(K_i)
  \right\},
\end{equation}
obtained by requiring that each function $v_h \in V_h$, when restricted to an
element $K_i$, lies in the corresponding local polynomial space
$\mathcal{Q}^p(K_i)$. We note that no continuity is enforced between the mesh
elements. We consider two neighboring elements, $K^-$ and $K^+$, that share a
face $e = \partial K^- \cap \partial K^+$. A function $v_h \in V_h$ is not
well-defined on $e$, and thus we define $v_h^-$ to be the trace of $v_h$ on $e$
from within $K^-$, and, analogously, $v_h^+$ to be the trace of $v_h$ on $e$
from within $K^+$. Similarly, $\bm n^-$ refers to a vector normal to $e$, facing
outward from $K^-$, and $\bm n^+ = -\bm n^-$ is the normal facing outward from
$K^+$. At this point, it will be useful to introduce the average $\{v_h\}$ and
jump $\llb v_h \rrb$ of a scalar function $v_h \in V_h$, which we define by
\begin{equation}
       \{ v_h \} = \frac{1}{2}\left( v_h^- + v_h^+ \right), \qquad
  \llb v_h \rrb  = v_h^- \bm n^- + v_h^+ \bm n^+.
\end{equation}
Similarly, for a vector-valued function $\bm\tau_h \in [V_h]^2$, we define
\begin{equation}
  \{ \bm\tau_h \} = \frac{1}{2}\left( \bm\tau_h^- + \bm\tau_h^+ \right), \qquad
  \llb \bm\tau_h \rrb  = \bm\tau_h^- \cdot \bm n^- + \bm\tau_h^+ \cdot \bm n^+.
\end{equation}
We remark that the jump of a scalar is a vector parallel to the normal
direction, and the jump of a vector is a scalar.

To obtain the DG discretization of equation~\ref{eq:heat-eqn-1}--\ref{eq:heat-eqn-3},
we transform the equation into a system of first order equations by introducing
the gradient $\bm\sigma = \nabla u$, thus obtaining the equivalent system of 
equations
\begin{align}
  \label{eq:heat-system-2}
  \bm\sigma &= \nabla u, \\
  \label{eq:heat-system-1}
  \frac{\partial u}{\partial t} &= \nabla \cdot \bm \sigma + f.
\end{align}
We look for an approximate solution $u_h \in V_h$, $\bm\sigma_h \in [V_h]^2$. To
obtain the variational form, we multiply equation \ref{eq:heat-system-1} by a
test function $v_h \in V_h$, and equation \ref{eq:heat-system-2} by a test
function $\bm\tau_h \in [V_h]^2$. We integrate the resulting equations over the
spatial domain $\Omega$. We then integrate the divergence and gradient terms by
parts element-by-element, giving rise to the following weak form,
\begin{gather} \label{eq:system-form-1}
    \int_{\Omega} \bm{\sigma}_h \cdot \bm{\tau}_h\,dx
    = -\int_{\Omega} u_h \nabla \cdot \bm{\tau}_h \,dx
    + \int_{\Gamma} \widehat{u}_h  \llb \bm{\tau}_h \rrb \, ds,
    \\ \label{eq:system-form-2}
    \int_\Omega \frac{\partial u_h}{\partial t} v_h \, dx =
    -\int_{\Omega} \bm{\sigma}_h \cdot \nabla v_h \, dx
    + \int_{\Gamma} \widehat{\bm{\sigma}}_h \cdot \llb v_h \rrb \, ds
    + \int_{\Omega} f v_h \, dx,
\end{gather}
where $\Gamma$ denotes all interior and exterior edges in the triangulation
$\mathcal{T}_h$, and $\widehat{u}_h$ and $\widehat{\bm\sigma}_h$ are
yet-to-be-defined \textit{numerical flux functions}. In this work, we consider
the local discontinuous Galerkin (LDG) method \cite{Cockburn1998a}, which 
proceeds by choosing the fluxes on interior edges by
\begin{align}
  \widehat{u}_h &= \{ u_h \} - \bm C_{12} \cdot \llb u_h \rrb,  \\
  \widehat{\bm\sigma}_h &= \{ \bm\sigma_h \} - C_{11} \llb u_h \rrb
    + \bm C_{12} \llb \bm\sigma_h\rrb,
\end{align}
for parameters $C_{11} \geq 0$ and $\bm C_{12}$. On the Dirichlet boundary, we
choose
\begin{align}
  \widehat{u}_h &= g_D, \\
  \label{eq:dirichlet-flux}
  \widehat{\bm\sigma}_h &= \bm\sigma_h^- - C_{11}(u_h^- - g_D)\bm n,
\end{align}
and on the Neumann boundary,
\begin{align}
  \widehat{u}_h &= u_h^-,\\
  \widehat{\bm\sigma}_h &= \bm g_N.
\end{align}
The parameter $C_{11}$ is necessary to stabilize the method, and can be thought
of as introducing an artificial viscosity \cite{Cockburn2001}. Of particular
interest to this work is the so-called minimal dissipation LDG method
\cite{Cockburn2007}, which allows for $C_{11}$ to be taken as identically zero
on all interior edges by a careful choice of the parameter $\bm C_{12}$.
Because the flux $\widehat{u}_h$ is independent of $\bm\sigma_h$, it is possible
to solve for $\sigma_h$ in an element-by-element fashion.

We choose a nodal basis for the space $V_h$, and write $\bm u$ to represent the
vector of degrees of freedom defining $u_h$. Likewise, $\bm\sigma$ represents
the vector of degrees of freedom defining $\bm\sigma_h$. We then rewrite
equations \ref{eq:system-form-1} and \ref{eq:system-form-2} as
\begin{align}
  M \bm \sigma &= G \bm u, \\
  M \bm u_t &= D \bm \sigma + E \bm u + M \bm f, 
\end{align}
or, equivalently,
\begin{equation}
  M \bm u_t = \left( D M^{-1} G + E \right) \bm u + M \bm f,
\end{equation}
where $M$ is the is standard mass matrix, and $D$ and $G$ are the divergence and
gradient operators defined by the respective bilinear forms. The mass matrix
times a vector-valued quantity is understood to be applied component-wise. We
note that a simple integration by parts shows that $G = -D^T$. The matrix $E$
corresponds to the stabilization terms caused by $C_{11} > 0$. This allows us to
define the Laplacian operator
\begin{equation}
\label{equ:DG_A_op}
  A = - D M^{-1} D^T + E,
\end{equation}
which is a symmetric, negative-definite linear operator. In the particular cases
of Neumann and periodic boundary conditions, we choose $C_{11}$ to be 
identically zero, and thus $E = 0$, so we can write $A = - D M^{-1} D^T$.

\subsection{Stochastic Driving Fields for DG} \label{sec:stochastic-fluxes}
\label{sec:stoch_driving_fields}
We now consider discretization of the stochastic diffusion
equation~\ref{equ:stoch_diffusion_main} where throughout we take diffusivity $D = 1$ and reference concentration $\bar{u} = 1$.  This corresponds to a stochastic version of equations \ref{eq:heat-eqn-1}--\ref{eq:heat-eqn-3}.  In equation~\ref{eq:heat-eqn-1}, we account for the fluctuations by taking the forcing term $f = \mb{g}$ to be a Gaussian stochastic field that is $\delta$-correlated in time with mean zero and with spatial covariance $\Lambda = \langle \bm f \bm f^T \rangle$.  The DG
discretization gives a linear stochastic dynamical system of the
form \ref{eq:governingEq_discr}, where the linear operator is given by 
\begin{eqnarray}
\label{equ:DG_L_op}
L = M^{-1} A
= -M^{-1} D M^{-1} D^T + M^{-1}E.
\end{eqnarray}
 The covariance structure of the fluctuations
is related to the linear operator $A$ and equilibrium covariance of the system
$C$ by the FDD framework through equation~\ref{equ:fluct_dissp_discr_space} by
\begin{equation}
  \Lambda = -L C - C^T L^T.
\end{equation}
At equilibrium, we have from equation~\ref{equ:stoch_diffusion_main} that $\langle u(\mb{x})u(\mb{y}) \rangle = \delta(\mb{x} - \mb{y})$.  In the discrete setting, the associated energy $\int_\Omega u_h^2\,dx = \bm{u}^T M \bm{u}$ motivates taking the discretized covariance to be $C = M^{-1}$ for the DG solution $\bm u$. From equation~\ref{equ:fluct_dissp_discr_space} and~\ref{equ:DG_L_op}, we use stochastic driving terms $\bm f$ with covariance
\begin{equation} \label{eq:f-covariance}
  \Lambda = \langle \bm f \bm f^T \rangle = 2 \left( M^{-1} D M^{-1} D^T M^{-1} - M^{-1} E M^{-1} \right).
\end{equation}
To obtain this expression, we have used that $M$ and $E$ are symmetric linear operators. By equation~\ref{equ:C_relate_G_spatial}, this prescribes for a given choice of mesh a covariance for the stochastic driving fields that takes into account the dissipative properties of the DG Laplacian operator to ensure at statistical steady-state the covariance $C$ is achieved.  

In order for this to be useful in practice, we must develop computational methods to efficiently generate the random forcing terms $\bm f$ with the prescribed covariance given by equation~\ref{eq:f-covariance}.  In principle, taking a Cholesky factorization $\Lambda = QQ^T$ would provide such random variates through $\bm f = Q \bs{\xi}$ with $\bs{\xi}$ a standard Gaussian random variable. However, this technique does not fully utilize the sparse block structure of the DG discretizations, and thus is often prohibitively expensive. Furthermore, such direct solvers can prove problematic to parallelize.

In the following sections, we shall develop more efficient methods for generating the stochastic variates, while also appropriately treating the boundary conditions, by utilizing the block structure of the matrices arising in DG.  Our methods scale as $\mathcal{O}\left(n_t (p+1)^{3d}\right)$, where the total number of degrees of freedom is given by $N = n_t (p+1)^d$,  $n_t$ is the number of mesh elements, $p$ the polynomial degree, and $d$ the spatial dimension.  We focus primarily on $h$-refinement where $p$ is held fixed and develop methods with linear computational complexity $\mathcal{O}(N)$ as $n_t \rightarrow \infty$.  Our methods have the additional advantage of being computed element-wise, avoiding communication between mesh elements, facilitating straightforward parallelization.

\subsubsection{Neumann and periodic boundary conditions with fluctuations} \label{sec:neumann}

We now consider the case of periodic or pure Neumann boundary conditions. In 
this case, we modify the finite element space to consist only of mean-zero 
functions,
\begin{equation}
  V_{h,0} = \left\{ v_h(\bm x) : \Omega \to \mathbb{R} : 
  v_h |_{K_i} \in \mathcal{P}^p(K_i), \text{ and }
  \int_\Omega v_h\,d\bm x = 0
  \right\},
\end{equation}
in order to ensure that the Laplacian operator $L$ is nonsingular. Other than
this modification, the method is as described in the preceding section. Using
the minimal dissipation LDG method, we can, in this case, set the LDG
stabilization parameter $C_{11}$ to be identically zero on all edges. As a
result, the stabilization matrix $E$ is also zero, and thus the covariance of
$\bm f$ is given by
\begin{equation}
  \Lambda = 2  M^{-1} D M^{-1} D^T M^{-1}.
\end{equation}
We recall that the DG mass matrix $M$ is given by
\begin{equation}
  M_{ij} = \int_\Omega \phi_i (\bm x) \phi_j(\bm x) \, d\bm x,
\end{equation}
where the functions $\phi_i$ form a basis for the finite element space $V_h$.
This matrix is symmetric positive-definite, and due to the discontinuous nature
of the basis functions $\phi_i$, has a natural block-diagonal structure, with blocks
corresponding to each element $K$ in the triangulation $\mathcal{T}_h$. Thus,
each block of the mass matrix is also symmetric positive-definite, allowing for
efficient block-wise computation of the Cholesky factorization $M ^{-1}= Q Q^T$
of the inverse mass matrix. 

These operations can be performed element-by-element and in parallel.  Given the block structure, the Cholesky factorization requires $\mathcal{O}(n_t (p+1)^{3d})$ operations.  The number of elements in the mesh is $n_t$, $p$ the polynomial order used, and $d$ the spatial dimension.  When performing $h$-refinement (fixing degree $p$), we have linear computational complexity with scaling $\mathcal{O}(N)$ where $N = n_t (p+1)^{3d} \sim C n_t$.  For $p$-refinement we see the complexity would be dimension $d$ dependent.

We can use this approach by defining the matrix $R$ by
\begin{equation}
  R = \sqrt{2} M^{-1} D Q,
\end{equation}
so that we obtain
\begin{equation}
  R R^T = 2 M^{-1} D Q Q^T D^T M^{-1} = \Lambda.
\end{equation}
We can generate variates using $\mb{f} = R\bs{\xi}$ where $\bm\xi$ is a Gaussian with $\delta$-correlation in time and spatial components with mean zero and variance one.  This has the covariance
\begin{equation}
 \lla \mb{f}\mb{f}^T \rra = \lla (R\bm\xi) (R\bm\xi)^T \rra =  R \lla \bm\xi \bm\xi^T \rra R^T =  R I R^T =\Lambda,
\end{equation}
thus giving a linear-time complexity method under $h$-refinement for generating
the random variates $\bm f$ prescribed by the FDD framework in
equation~\ref{equ:fluct_dissp_discr_space} and~\ref{eq:f-covariance}.

\subsubsection{Dirichlet boundary conditions with fluctuations} 
\label{sec:dirichlet}

We now consider the case of Dirichlet boundary conditions, which are enforced
by penalizing the difference between the approximate solution $u_h$ and the 
prescribed boundary value $g_D$ by means of a penalty parameter $C_{11}$, as in
equation \ref{eq:dirichlet-flux}. In order for the problem to be well-posed,
we must choose $C_{11} > 0$ on all edges in $\partial \Omega$. As before, we 
set $C_{11} = 0$ on all interior edges. This results in a nonzero penalty matrix
$E$, whose entries are given by
\begin{equation}
  E_{ij} = - \int_{\partial\Omega} C_{11} \phi_i(\bm{x}) \phi_j(\bm{x}) \, ds.
\end{equation}
This matrix is symmetric and negative-semidefinite, with an element-wise block 
diagonal structure. The nonzero blocks of $E$ correspond to elements of the 
mesh with at least one edge on the domain boundary.

We remark that the discontinuous Galerkin method enforces the Dirichlet boundary
conditions weakly through the penalization procedure described above. This 
formulation can permit fluctuations at the domain boundary. In this case, we
consider the target covariance to be given by $C = M^{-1}$. The FDD framework based on fluctuation-dissipation balance in equation~\ref{equ:fluct_dissp_discr_space} gives the prescribed covariance in equation~\ref{eq:f-covariance}. To efficiently sample a random variable with this covariance, we require the eigendecomposition of the penalty matrix,
\begin{equation}
  E = V D V^T.
\end{equation}
This factorization can be computed efficiently in linear-time complexity under $h$-refinement by taking advantage of the block structure of the matrix $E$. We then let $\bm\xi_1$ and $\bm\xi_2$ be
independent, identically distributed Gaussian random variables with mean zero
and variance one. We define the matrix $R_1$ by
\begin{equation}
  R_1 = \sqrt{2} M^{-1} D Q,
\end{equation}
where as in the Neumann case, $Q$ is the Cholesky factor of the inverse mass 
matrix. We also define $R_2$ by
\begin{equation}
  R_2 = \sqrt{2} M^{-1} V \sqrt{-D}.
\end{equation}
Then, by setting $\bm f = R_1 \bm\xi_1 + R_2 \bm\xi_2$, we obtain
\begin{equation}
  \lla \bm f \bm f^T \rra = 2 M^{-1} D M^{-1} D^T M^{-1}
    M^{-1} + 2 V (-D) V^T M^{-1} = \Lambda.
\end{equation}

This formulation leads to fluctuations on the domain boundary, 
which may be considered undesirable. If we wish not to introduce fluctuations
at the boundary, we may modify the target covariance as follows. We recall that
our basis functions $\phi_i$ are defined to be nodal interpolants. We refer to
an index $i$ as a \textit{boundary index} if the nodal interpolation point 
corresponding to basis function $\phi_i$ lies on the domain boundary.
We then define the modified target covariance matrix $\widetilde{C}$ by
\begin{equation}
  \widetilde{C}_{ij} = \begin{cases}
    0 &  \text{if either $i$ or $j$ is a boundary index,} \\
    (M^{-1})_{ij} & \text{otherwise}.
  \end{cases}
\end{equation}
In this case, the fluctuation covariance is given by
\begin{equation}
  \Lambda = -L\widetilde{C} - \widetilde{C}^T L^T.
\end{equation}
Here, we note that $-L\widetilde{C}$ is no longer symmetric. In order to ensure symmetry,
we modify the Laplacian operator as follows. Rather than weakly enforce the 
Dirichlet conditions using a penalty term, we strongly enforce these conditions
by fixing the degrees of freedom at the boundary to be equal to their specified
Dirichlet values, $g_D$. The resulting linear operator $\widetilde{L}$ has the
form
\begin{equation}
  \widetilde{L} = -\widetilde{I} M^{-1} D M^{-1} D^T \widetilde{I},
\end{equation}
where the matrix $\widetilde{I}$ is a diagonal matrix defined by
\begin{equation}
  \widetilde{I}_{ii} = \begin{cases}
    0 &  \text{if $i$ is a boundary index,} \\
    1 & \text{otherwise}.
  \end{cases}
\end{equation}
We note that $\widetilde{I}E = 0$, and thus the boundary penalty matrix is not needed in
this formulation. Additionally, we have that $\widetilde{I} M^{-1} = \widetilde{C}$ and $\widetilde{I}
\widetilde{C} = \widetilde{C}$. Thus, we have
\begin{equation}
  \widetilde{\Lambda}
    = -\widetilde{L}\widetilde{C} - \widetilde{C}^T\widetilde{L}^T
    = \widetilde{I} M^{-1} D M^{-1} D^T \widetilde{I}\widetilde{C}
    + \widetilde{C}^T\widetilde{I} M^{-1} D M^{-1} D^T \widetilde{I}
    = 2 \widetilde{C} D M^{-1} D^T \widetilde{C}.
\end{equation}
Thus, in this case, we can define the matrix $\widetilde{R}$ by
\begin{equation}
  \widetilde{R} = \sqrt{2} \widetilde{C} D Q,
\end{equation}
where as before $Q$ is the Cholesky factor of the inverse mass matrix. We 
define the fluctuation forcing term by $\bm f = \widetilde{R} \bm\xi$, and thus
obtain
\begin{equation}
  \lla \bm f \bm f^T \rra =   \widetilde{R}  \widetilde{R}^T = \Lambda.
\end{equation}
Given the block structure each required factorization operation only costs linear-time computational complexity under $h$-refinement.  In this way, we can efficiently generate in practice the needed stochastic driving terms with the prescribed covariance structure given by the FDD framework in equation~\ref{equ:fluct_dissp_discr_space} and~\ref{eq:f-covariance}.

\subsubsection{Temporal discretization}
We use an Euler-Maruyama time discretization \cite{Kloeden1992}, which
results in the fully-discrete system given by
\begin{equation}
  \bm u^{n+1} = \bm u^n + L \bm u^n \Delta t + \bm f^n,
\end{equation}
with time step $\Delta t$. The term $\bm f^n$ is used to denote a Gaussian 
random variable with mean zero and covariance
\begin{equation}
  \lla (\bm f^m)(\bm f^n)^T \rra = \Lambda\Delta t \delta_{mn}.
\end{equation}
This random variable is generated at each time step by sampling the random
fluctuation forcing term $\bm f^n \sim \sqrt{\Delta{t}}\bm f$ according to the desired boundary conditions,
using the methodology described in Sections \ref{sec:neumann} and
\ref{sec:dirichlet}. 

We mention that while we could in principle use the relation in equation~\ref{equ:fluct_dissp_discr_time} to make further corrections to the covariance of the stochastic driving terms $\Lambda$.  In practice, the time steps $\Delta{t}$ we take in our numerical calculations are sufficiently small that these higher-order corrections do not play a significant role. Therefore, for the sake of simplicity, we do not incorporate these terms into our discretization in this work.  We also remark that other stochastic time-step integrators could also be developed with relations derived like in equation~\ref{equ:fluct_dissp_discr_time} to make further corrections to the covariance of the stochastic driving terms.

\section{Results}
\label{sec:results}

We investigate how these methods perform in practice by considering problems involving different geometries and boundary conditions.  We show how the stochastic discretizations perform when accounting for the Neumann and Dirichlet boundary conditions.  We also consider how the stochastic methods perform on different types of meshes ranging from structured to unstructured, and when transforming from straight-sided meshes to curved geometries.  We investigate how the methods behave using both low-order and high-order polynomial spaces.  These investigations illustrate some of the key features of our stochastic DG methods.

We compare the results of our fluctuation-dissipation based discretizations for our stochastic DG methods described in Section~\ref{sec:stochastic-fluxes} with those that would be obtained from the intuitive approach of directly using random fluxes.  The random fluxes are given by $F
= M^{-1} D \bm \xi$, where $\bm \xi$ is the vector whose entries are given by
independent, identically-distributed Gaussian random variables, and represents a
random function $\xi_h \in V_h$ when expanded in terms of the chosen basis of the function space $V_h$.  To characterize how stochastic features of the system are captured by the stochastic numerical methods we investigate the covariance structure of fluctuations using Monte-Carlo sampling with estimator
$C \approx \frac{1}{N} \sum_{k=1}^N (\bm u^k) (\bm u^k)^T$.
We take only samples with $t \geq t_0$ for $t_0$ typically covering about $10\%$ of our samples to exclude initial transients in the stochastic dynamics.  When considering the stochastic systems we assume throughout ergodicity.

We remark that in our methods if, instead of using a nodal basis, we were to use a basis consisting of orthogonal polynomials on each element, then the resulting mass matrix would be the identity matrix, and the random fluxes described above would agree with those obtained using the fluctuation-dissipation principle.  This would hold except in the case of weakly-enforced Dirichlet conditions, in which case the
boundary penalization term must also be taken into account. In order to ensure that the
units are consistent between our comparisons, we multiply the random fluxes by a geometric factor of $(p+1)^2 / h$, which represents the spatial resolution of the DG method.

We compare the results of our fluctuation-dissipation based SDGM and what would be obtained for random fluxes in the cases of periodic, Neumann, and Dirichlet boundary conditions.  These present a strong test of how methods capture inherent statistical structure in these stochastic systems.

\subsection{Periodic Boundary Conditions}
\label{sec:results_periodic}

\begin{figure}[H]
  \centering
  \raisebox{-0.5\height}{\includegraphics[width=0.8\columnwidth]{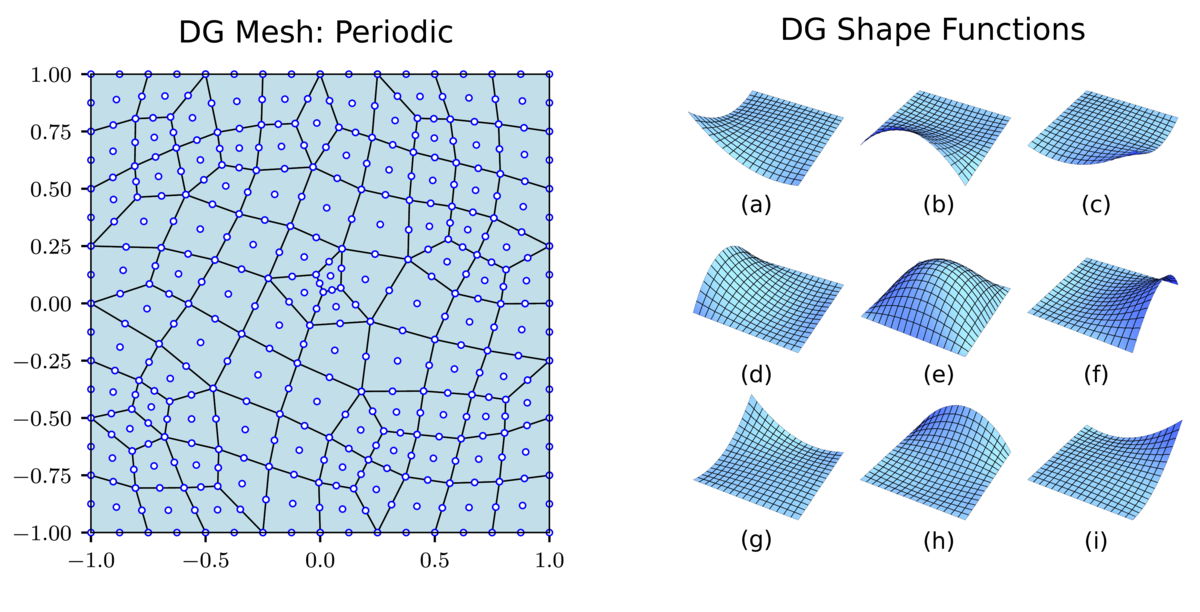}}
  \caption{Left: Unstructured mesh of $[-1,1] \times [-1,1]$ with periodic
  boundary conditions along edges. DG nodes corresponding to $p=2$ are shown in blue. Right (a)--(i): 
  Biquadratic DG shape functions for Gauss-Lobatto nodal basis on each element.}
  \label{fig:periodic-mesh}
\end{figure}

\begin{figure}[t]
  \centering
  \includegraphics[width=0.8\columnwidth]{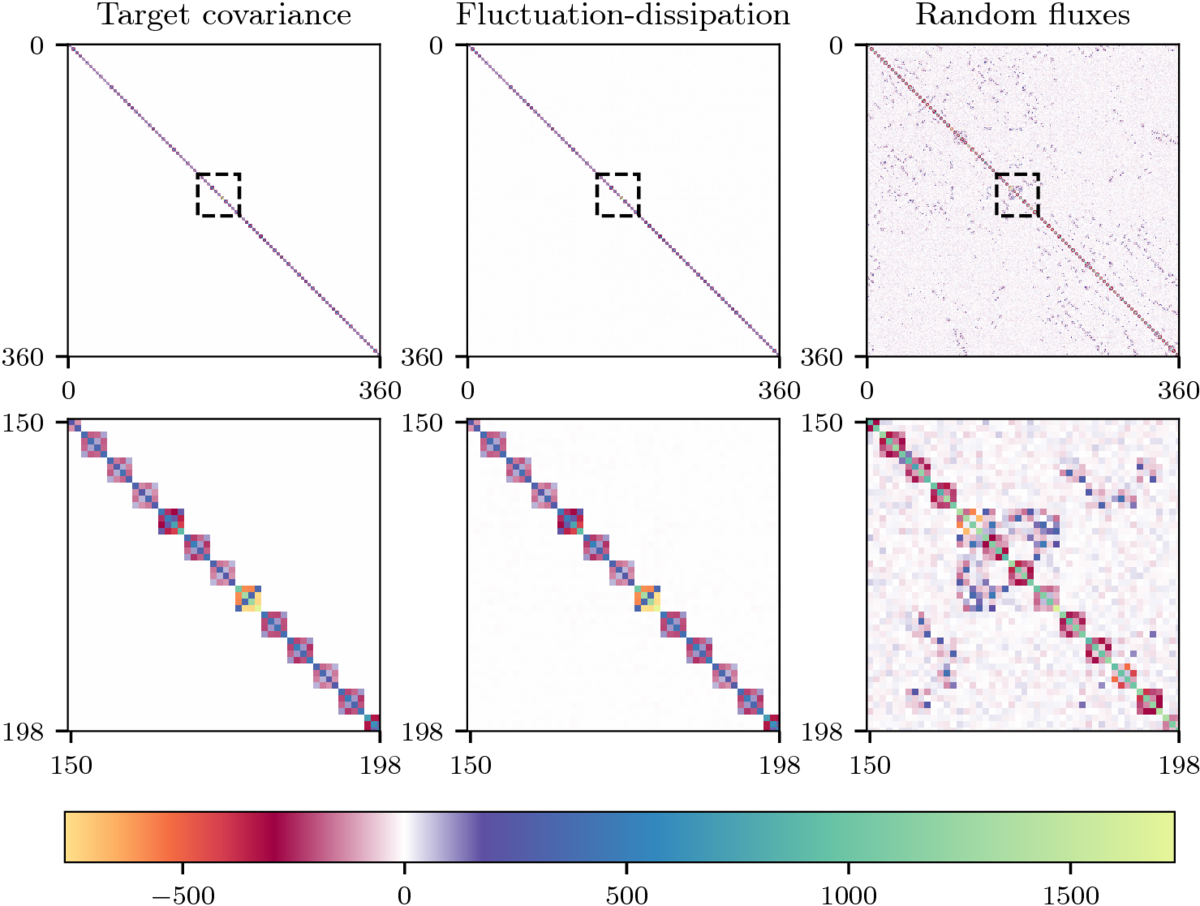}
  \caption{Comparison of target covariance (left) with empirical covariance 
           matrices obtained using the fluctuation-dissipation principle
           (center) and random fluxes (right) with periodic boundary conditions.
           A zoom-in is shown on the bottom row for comparison.}
  \label{fig:periodic-matrices}
\end{figure}

We consider the case of periodic boundary conditions for the domain having a simple geometry given by $\Omega = [-1,1] \times
[-1,1]$.  We enforce periodic boundary conditions on the edges using symmetry of the nodes.  Relative to many other discretization methods, an advantage of the DG method is the ability to use general, unstructured meshes.  To illustrate this property of the method, we use the mesh shown in Figure~\ref{fig:periodic-mesh}.  We choose the finite element function space to consist
of $p=1$ piecewise bilinear polynomials in each element, thus resulting in a
method whose spatial convergence rate is $\mathcal{O}(h^2)$. The function space
$V_h$ is restricted to the space of mean-zero functions on $\Omega$ to ensure
non-singularity of the Laplacian operator.

\begin{figure}[H]
  \centering
  \includegraphics[width=0.7\columnwidth]{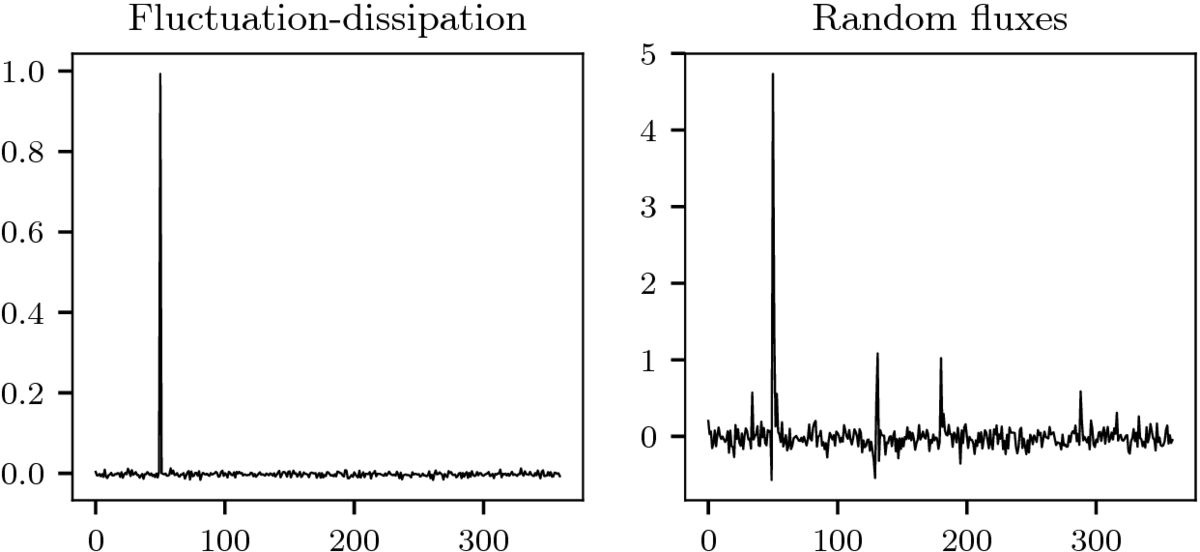}
  \caption{Correlations of 50th nodal degree of freedom, periodic
           test case. Fluctuation-dissipation forcing term (left) shows clean
           $\delta$-correlation, random fluxes shown on right.}
  \label{fig:periodic-row}
\end{figure}

\begin{figure}[H]
  \centering
  \includegraphics[width=0.7\columnwidth]{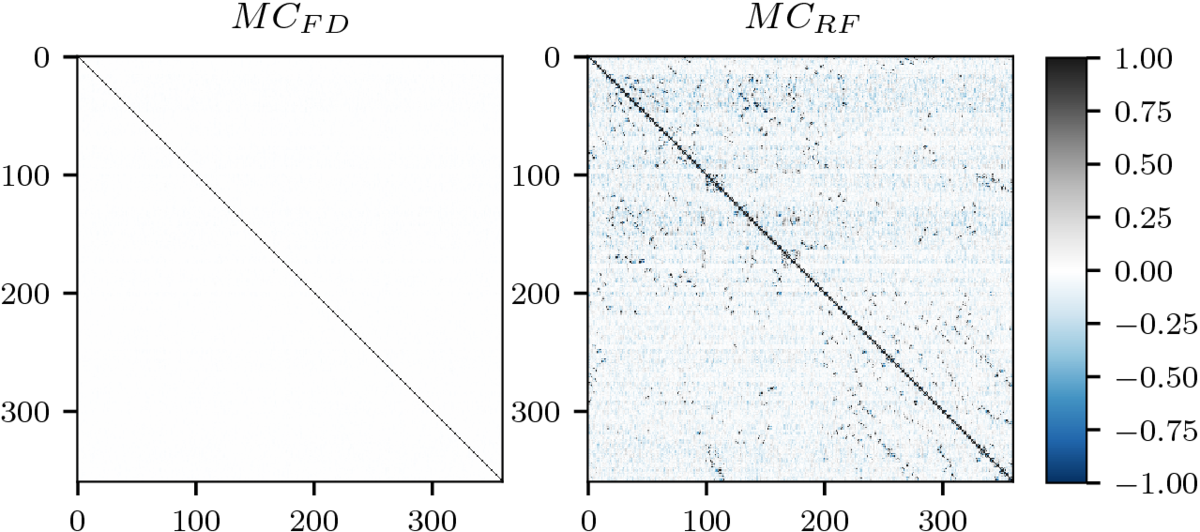}
  \caption{Comparison of $M C_{FD}$ (left) with $M C_{RF}$ (right) for 
           periodic boundary conditions. Both matrices are supposed to
           approximate the identity matrix. Significant spurious long-range
           correlations are seen in $MC_{RF}$.}
  \label{fig:periodic-matrices-2}
\end{figure}

We choose the time step to be given by $\Delta t = 10^{-5}$, and approximate the
covariance of the solution using 
$5\times 10^6$
Monte Carlo samples. We let
$C_{FD}$ denote the covariance of the solution obtained using forcing terms
determined according to the fluctuation-dissipation principle, and let $C_{RF}$
denote the covariance of the solution obtained using a random flux forcing
function. We compare these two empirical covariance matrices with the target
covariance, which is given by the inverse mass matrix, $C = M^{-1}$. The
relative $L^2$ error $\| C_{FD} - C \|_{L^2} / \| C \|_{L^2}$ for the
fluctuation-dissipation case was 3.5\%. We compare the resulting matrices in
Figure~\ref{fig:periodic-matrices}. We notice that the matrix obtained using the
fluctuation-dissipation principle displays excellent agreement with the target
covariance matrix, with close to no correlations in between elements. On the
other hand, the covariance matrix obtained using random fluxes does not
accurately reproduce the target covariance, and displays significant spurious
long-range correlations. For ease of comparison, we note that $M C = I$, and
thus we expect $M C_{FD}$ and $M C_{RF}$ to well-approximate the identity matrix
as well. We show both these matrices in Figure~\ref{fig:periodic-matrices-2}. We
note that $M C_{RF}$, resulting from the use of random fluxes, does not display
good agreement with the identity matrix, illustrating the spurious correlations
that arise from this strategy.

\subsection{Neumann Boundary Conditions}
\label{sec:results_neumann}

\begin{figure}[H]
  \centering
  \includegraphics[width=0.4\columnwidth]{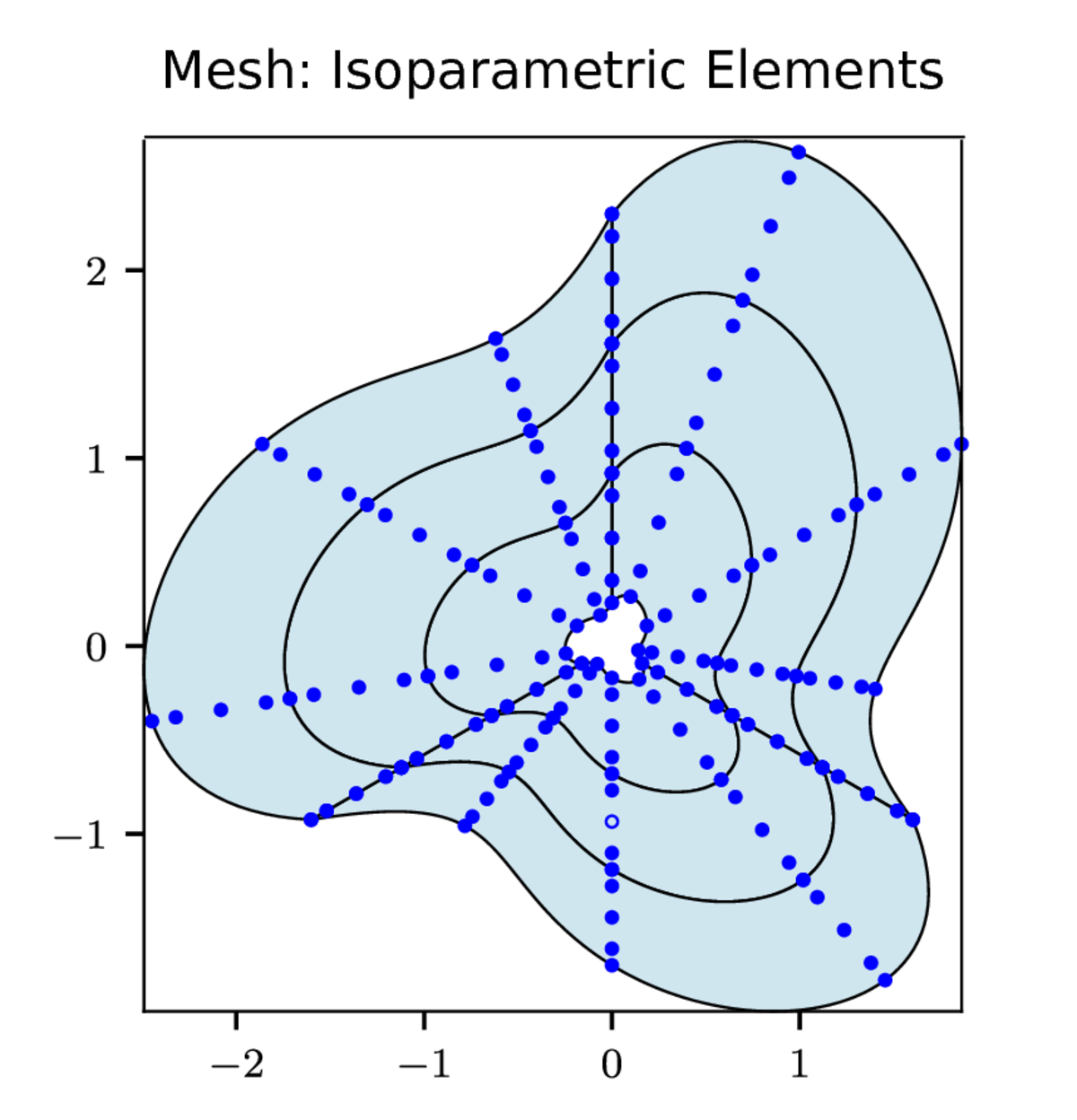}
  \caption{High-order mesh of warped annulus, with isoparameteric curved 
  elements. Blue circles indicate $p=4$ DG nodes.}
  \label{fig:annulus-mesh}
\end{figure}

Now, we extend the numerical tests beyond the case of periodic boundary
conditions. We consider a curved geometry that is given by a warped annulus. We
use a $p=4$ high-order polynomial basis. The mesh used for this problem is shown
in Figure~\ref{fig:annulus-mesh}. This examples illustrates the use of Neumann
boundary conditions, curved geometries via an isoparametric mapping, and the use
of high-order polynomial spaces. As in the case of periodic boundary conditions,
we restrict the finite element space to mean-zero functions, so that the Laplace
operator is nonsingular. We note again that the minimal dissipation LDG method
allows for the stabilization parameter $C_{11}$ to be identically zero on all 
edges in the mesh. Homogeneous Neumann conditions are enforced at the domain
boundary.

\begin{figure}[H]
  \centering
  \includegraphics[width=0.7\columnwidth]{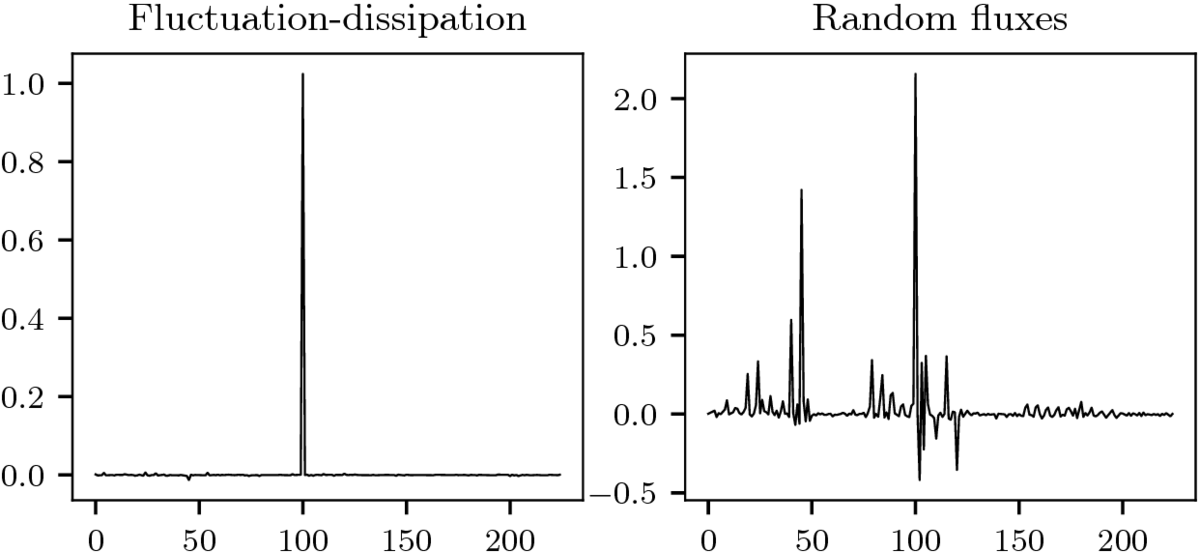}
  \caption{Correlations of 100th nodal degree of freedom, homogeneous Neumann
           test case. Fluctuation-dissipation forcing term (left) shows clean
           $\delta$-correlation, random fluxes shown on right. }
  \label{fig:neumann-row}
\end{figure}

\begin{figure}[H]
  \centering
  \includegraphics[width=0.8\columnwidth]{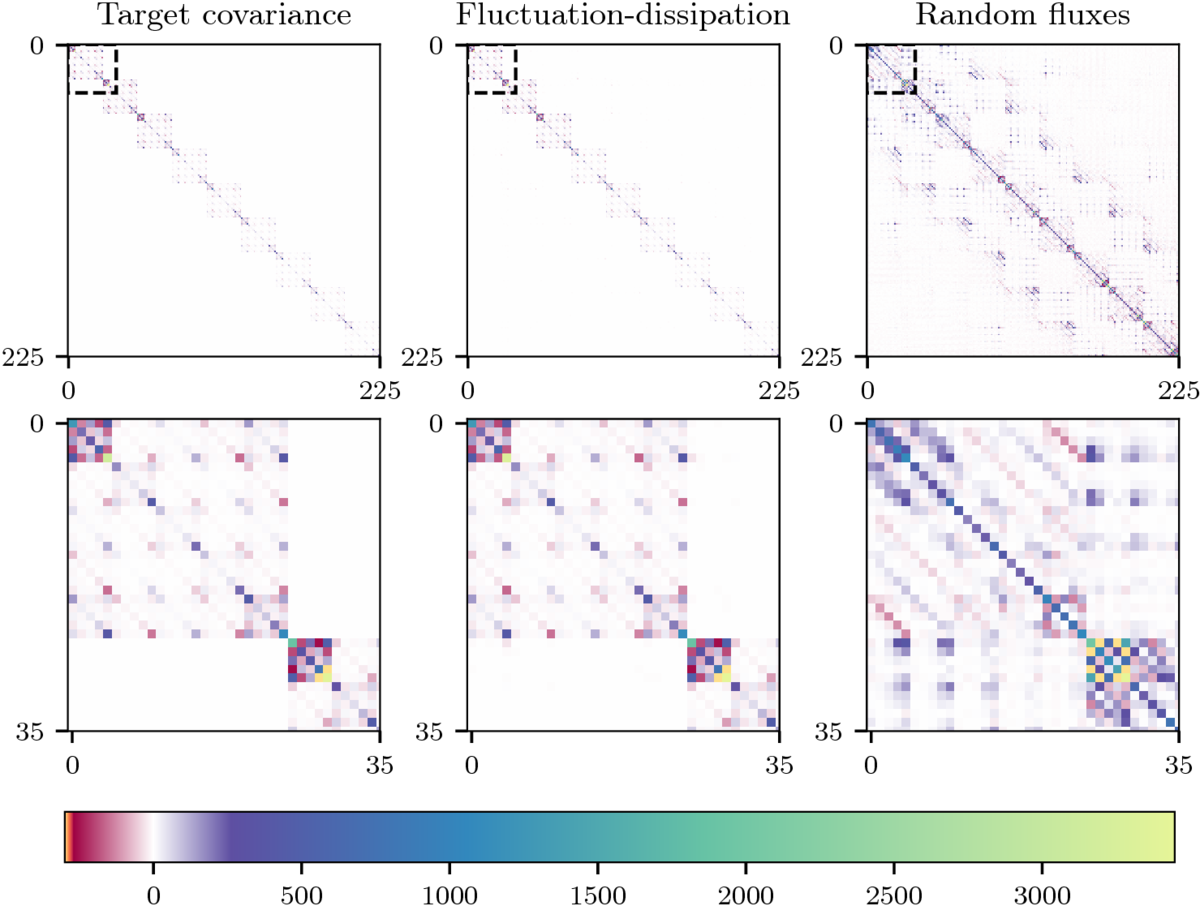}
  \caption{Comparison of target covariance (left) with empirical covariance 
           matrices obtained using the fluctuation-dissipation principle
           (center) and random fluxes (right) with Neumann boundary conditions.
           A zoom-in is shown on the bottom row for comparison.}
  \label{fig:neumann-matrices}
\end{figure}

$•$

We choose the time step to be $\Delta t = 10^{-5}$ and perform the same
comparisons as in the previous case. We estimate the covariance matrix using
$5\times 10^6$
samples. The relative $L^2$ error for the fluctuation-dissipation
case was 4.8\%. Comparisons of the matrices are shown in Figure~\ref{fig:neumann-matrices}. Additionally, in Figure~\ref{fig:neumann-row},
display a single row of the the matrices $M C_{FD}$ and $M C_{RF}$, which
represents the correlation between a specified nodal degree of freedom with all
other nodal degrees of freedom. This figure illustrates the clean
$\delta$-correlation that arises from using the fluctuation-dissipation forcing
terms. On the other hand, we see that using random flux forcing terms does not
result in the desired $\delta$-correlation.

\subsection{Dirichlet Boundary Conditions}
\label{sec:results_dirichlet}

\begin{figure}[H]
  \centering
  \includegraphics[width=0.8\columnwidth]{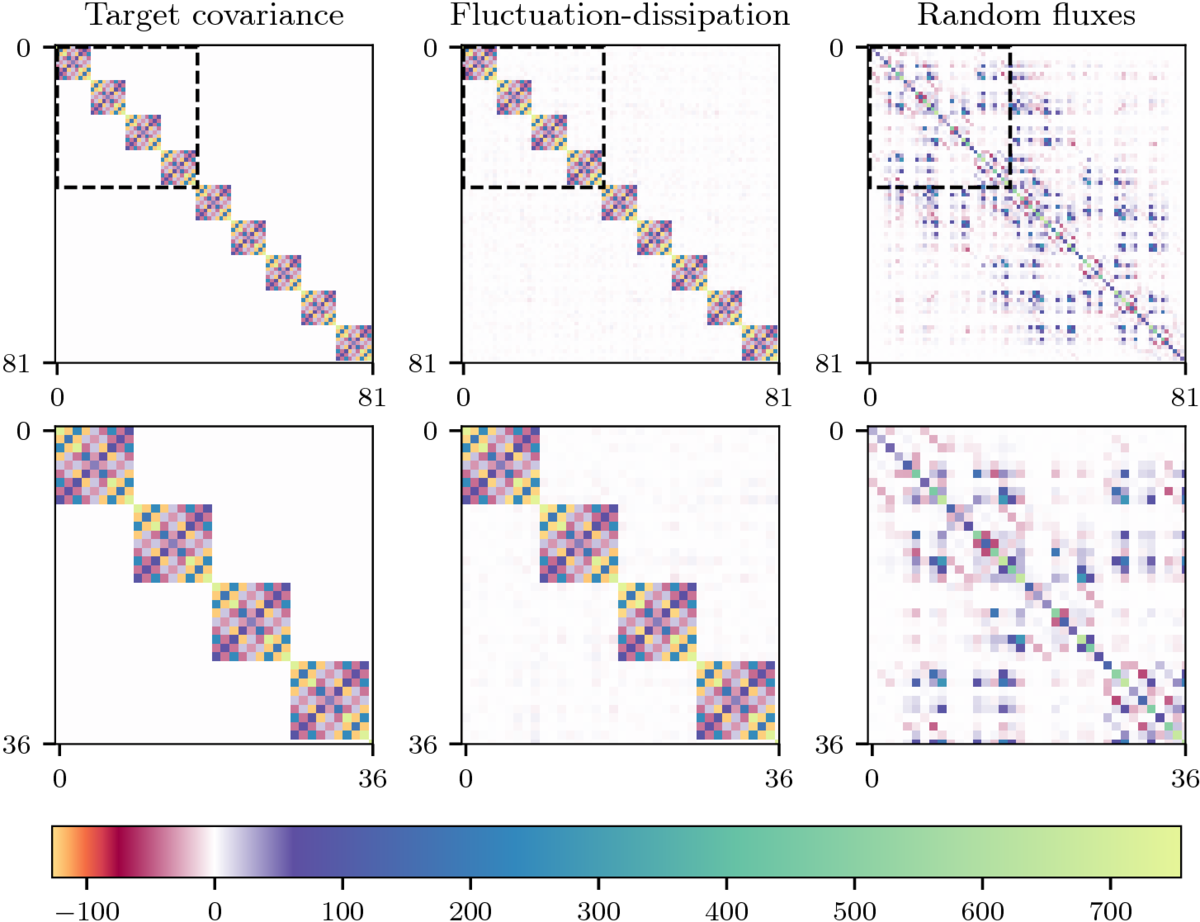}
  \caption{Comparison of target covariance (left) with empirical covariance 
           matrices obtained using the fluctuation-dissipation principle
           (center) and random fluxes (right) with weakly-imposed Dirichlet
           conditions. A zoom-in is shown on the bottom row for comparison.}
  \label{fig:dirichlet-weak-matrices}
  \end{figure}

For a final test case, we consider $\Omega = [0,1] \times [0,1]$, with
homogeneous Dirichlet conditions imposed on $\partial \Omega$. We discretize the
geometry using a regular Cartesian grid, and use $p=2$ biquadratic polynomials.
First we consider allowing fluctuations on the domain boundary, which we achieve
by weakly enforcing the Dirichlet conditions using the penalty parameter, which
we choose to be strictly positive on all edges lying on $\partial \Omega$, and
equal to zero on all interior edges. The target covariance for this case is, as
in the preceding cases, $C = M^{-1}$. This case requires treatment of the
stabilization term $E$, as described in Section~\ref{sec:dirichlet}. Using
$5\times 10^6$
Monte Carlo samples, we find the relative error between $C$ and
$C_{FD}$ in the $L^2$ matrix norm is 3.5\%. A comparison of the empirical
covariance matrices for this case is shown in Figure~\ref{fig:dirichlet-weak-matrices}

\begin{figure}[H]
  \centering
  \includegraphics[width=0.8\columnwidth]{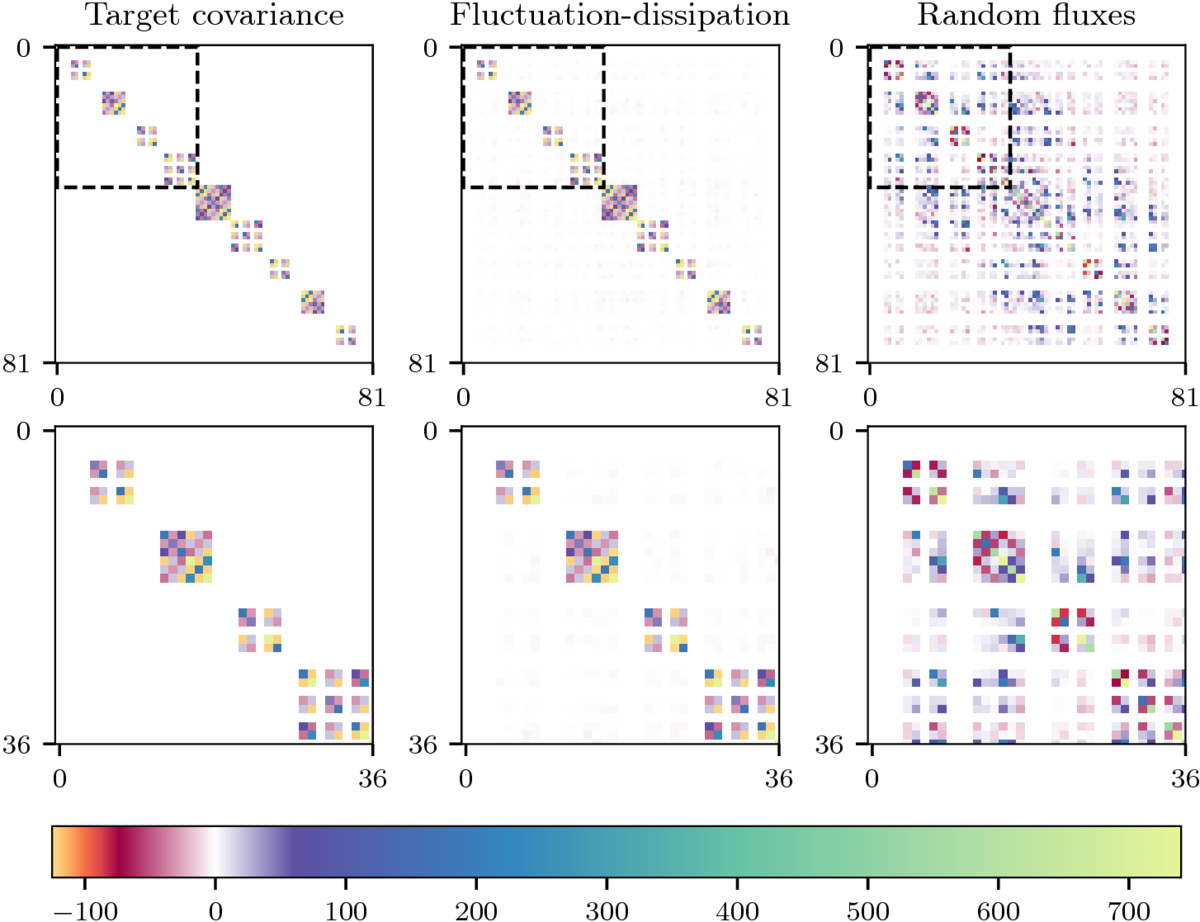}
  \caption{Comparison of target covariance (left) with empirical covariance 
           matrices obtained using the fluctuation-dissipation principle
           (center) and random fluxes (right) with strong Dirichlet conditions.
           A zoom-in is shown on the bottom row for comparison.}
  \label{fig:dirichlet-strong-matrices}
\end{figure}

We additionally consider the case where no fluctuations are permitted on the
domain boundary. The target covariance in this case is given by $\widetilde{C}$,
which is obtained from $C = M^{-1}$ by setting $\widetilde{C}_{ij} = 0$ whenever
either $i$ or $j$ is the index of a boundary node. To ensure symmetry of the 
resulting covariance matrix $\widetilde{\Lambda}$, we enforce the homogeneous 
Dirichlet conditions strongly rather than through a penalty term. We compare 
the resulting covariance matrices in Figure~\ref{fig:dirichlet-strong-matrices}.

\section{Conclusions}
\label{sec:conclusions}
We have introduced a general approach referred to as Fluctuation Dissipation Discretizations (FDD) for approximating SPDEs using fluctuation-dissipation balance.  Our FDD framework provides for a given numerical method's spatial and temporal discretizations a prescribed way to discretize the stochastic driving terms in order to take into account the dissipative properties of the discrete numerical operators and how this affects propagation of fluctuations in the system.  Using these ideas, we have developed general Stochastic Discontinuous Galerkin Methods (SDGM) for approximating solutions of parabolic stochastic partial differential equations.  We developed  practical methods with linear-time computational complexity under $h$-refinement for generating the needed random variates with the prescribed covariance structure for the DG discretizations on general geometries.  We have also introduced methods for handling periodic, Neumann, or Dirichlet boundary conditions.  We expect many of the ideas behind our FDD framework to provide helpful guidelines in the further development of SDGMs and other numerical discretizations for approximating stochastic partial differential equations.  

\section{Acknowledgments}
We acknowledge support to P.J.A. and N.T. from research grant DOE ASCR CM4 DE-SC0009254.  N.T. was also supported by the NSF-MSPRF program, and P.J.A. was also supported by NSF Grant DMS - 1616353.  W.P.\ was supported by the Department of Defense through the National Defense Science and Engineering Graduate Fellowship Program and by the Natural Sciences and Engineering Research Council of Canada. 

\bibliographystyle{plain}
\bibliography{paperDatabase}{}

\appendix

\end{document}